\documentclass[reqno,centertags, 12pt]{amsart}

\usepackage{amssymb,amsmath,amsfonts,amssymb} 
-\marginparwidth \oddsidemargin -4mm \evensidemargin
\oddsidemargin
\usepackage{latexsym}
\advance\hoffset by 5mm

\def\@abssec#1{\vspace{.05in}\footnotesize \parindent .2in
{\bf #1. }\ignorespaces}
\newtheorem{theorem}{Theorem}[section]
\newtheorem{thm}[theorem]{Theorem}

\newtheorem{corollary}[theorem]{Corollary}

\def \Rm {\mathbb R}

\allowdisplaybreaks \numberwithin{equation}{section}

\title{Spectral Properties of Schr\"odinger Operators with Decaying Potentials}
\author{Sergey A. Denisov and Alexander Kiselev}
\thanks{Department of
Mathematics, University of Wisconsin, Madison, WI 53706, USA;
e-mail: denissov@math.wisc.edu, kiselev@math.wisc.edu}

\begin{document}

\begin{abstract}
We review recent advances in the spectral theory of Schr\"odinger
operators with decaying potentials. The area has seen spectacular
progress in the past few years, stimulated by several conjectures
stated by Barry Simon starting at the 1994 International Congress
on Mathematical Physics in Paris. The one-dimensional picture is
now fairly complete, and provides many striking spectral examples.
The multidimensional picture is still far from clear and may
require deep original ideas for further progress. It might hold
the keys for better understanding of a wide range of spectral and
dynamical phenomena for Schr\"odinger operators in higher
dimensions.
\end{abstract}

\maketitle

\section{Introduction}

The Schr\"odinger operators with decaying potentials are used to
study the behavior of a charged particle in a local electric
field. The operator is defined by
\begin{equation}\label{SO}
H_V = -\Delta +V(x)
\end{equation}
on $L^2(\Rm^d);$ in one dimension it is common to consider the
operator on half-axis with some self-adjoint boundary condition at
zero. The spectral and dynamical effects that we are interested in
are those depending on the rate of decay of the potential rather
than its singularities, so we will often freely assume that $V$ is
bounded. If decay of the potential is sufficiently fast (short
range), one expects scattering motion. The corresponding results
have been rigorously proved by Weidmann \cite{Weid} in one
dimension (where short-range case means that the potential is in
$L^1(\Rm)$) and by Agmon \cite{Ag} in higher dimensions (where the
natural short range class is defined by $|V(x)| \leq
C(1+|x|)^{-1-\epsilon}$). These results established pure absolute
continuity of the spectrum on the positive semi-axis and
asymptotic completeness of the wave operators. There has been a
significant amount of work on longer range potentials with
additional symbol-like conditions (see e.g. \cite{RS3} for further
references), or oscillating potentials of specific Wigner-von
Neumann type structure (see e.g. \cite{BA,D1,White} for further
references). However until 1990s there has been very limited
progress on understanding slowly decaying potentials with no
additional assumptions on behavior of derivatives. The short range
or classical WKB methods did not seem to apply in this case, and
the possible spectral properties remained a mystery. The
celebrated Wigner-von Neumann example \cite{WvN} provides a
potential $V(x)$ behaving like $8\sin (2x) /x +O(x^{-2})$ as $x
\rightarrow \infty$ and leading to an imbedded eigenvalue $E=1,$
thus showing that surprising things can happen once potential is
not short range. On the other hand, the work of Kotani and
Ushiroya \cite{KU} implied that for potentials decaying at power
rate slower than $x^{-\alpha},$ $\alpha \leq 1/2,$ the spectrum
may become purely singular, and thus the scattering picture may be
completely destroyed. There was a clear gap in the decay rates
where very little information on the possible spectral properties
was available. In the recent years, there has been a significant
progress in the area, largely stimulated by Barry Simon's research
and ideas. At the ICMP in Paris in 1994, Simon posed a problem of
understanding the spectral properties of Schr\"odinger operators
with potentials satisfying $|V(x)| \leq C(1+|x|)^{-\alpha},$ $1
>\alpha > 1/2.$ Later, at the 2000 ICMP in London \cite{Simrev},
he compiled a list of fifteen problems in Schr\"odinger operators
``for the twenty first century". Two of the problems on the list
concern long range potentials.

While there remain many open questions, the recent effort to
improve understanding of the long range potentials led to many
high quality mathematical works. Fruitful new links between
spectral theory of Schr\"odinger operators and orthogonal
polynomials as well as Fourier analysis have been discovered and
exploited. Surprising examples of intricate spectral properties
have been produced. Advances have been made towards better
understanding of effects possible in higher dimensions. In this
review, we try to survey recent results in this vital area, as
well as underline most active current directions and questions of
interest. In the second section we discuss the one-dimensional
case, where the picture is much more detailed and complete. The
third section is devoted to a number of interesting spectral and
dynamical examples, typically one-dimensional, but easily
extendable to any dimension by spherically symmetric construction.
We briefly mention certain relations to Dirac operators, Jacobi
matrices and polynomials orthogonal on the unit circle (OPUC) in
the fourth section. In the last section, we consider the higher
dimensional case, where the main question - known as Barry Simon's
conjecture - is still open and is at the focus of current
research.

Most of this work is a compilation and review of known results.
There are, however, three nuggets that are new to the best of our
knowledge. In Section 3, we provide a new proof of
Theorem~\ref{dpis}, the construction of an example with a dense set
of imbedded eigenvalues. In Section 5, Theorem~\ref{bounsol} and
Theorem~\ref{improve} are new.

\section{The One-Dimensional Case}

The initial progress in understanding slowly decaying potentials
started from particular cases, such as random and sparse. Both of
these classes have been treated in a single framework by Barry
Simon in a joint paper with Last and AK \cite{KLS}; random
decaying potentials in the discrete setting have been pioneered by
Simon in a joint work with Delyon and Souillard \cite{DSS}. Let
$H_V$ be a half line Schr\"odinger operator, and fix some boundary
condition at the origin. Let us call $V(x)$ a Pearson potential if
$V(x) = \sum_n a_n W(x-x_n),$ where $a_n \rightarrow 0,$ $x_n
/x_{n-1} \rightarrow \infty$ as $n \rightarrow \infty,$ and $W(x)
\in C_0^\infty(\Rm).$
\begin{thm}\label{sparse1}
Let $V(x)$ be a Pearson potential. If $\sum_n a_n^2 < \infty,$ the
spectrum of $H_V$ on $(0,\infty)$ is purely absolutely continuous.
If $\sum_n a_n^2 = \infty,$ the spectrum of $H_V$ on $(0,\infty)$
is purely singular continuous.
\end{thm}
This result \cite{KLS} generalizes the original work of David
Pearson \cite{Pear}, who essentially proved Theorem~\ref{sparse1}
under assumption that $x_n$ grow sufficiently fast (with no
explicit estimate). See also \cite{Rem1} for related results. In a
sense, Pearson's theorem was the first indication of a clear
spectral transition at $p=2$ when potential is viewed in $L^p$
scale. A similar picture is true for the random potentials. Let
$V(x) = n^{-\alpha}a_n(\omega)W(x-n),$ where $W \in C_0^\infty
(0,1)$ and $a_n(\omega)$ are random i.i.d variables with mean zero
and compactly supported probability density function.
\begin{thm}\label{random1}
If $\alpha > 1/2,$ then the spectrum of $H_V$ on the positive half
axis is purely absolutely continuous with probability one. If
$\alpha < 1/2,$ the spectrum on $(0,\infty)$ is pure point with
probability one. If $\alpha =1/2,$ the spectrum is a mixture of
pure point and singular continuous spectrum with probability one.
\end{thm}
Please see \cite{KLS} for more details in the $\alpha =1/2$ case,
as well as for the proof of a more general theorem. The first
result of the type of Theorem~\ref{random1} has been due to
Delyon, Simon and Souillard \cite{DSS}, who handled the discrete
case. In the continuous setting, Kotani and Ushiroya \cite{KU}
proved a version of Theorem~\ref{random1} for a slightly different
model.

Theorems~\ref{sparse1}, \ref{random1} show the transition which is
reminiscent of some classical results on almost everywhere
convergence and divergence of Fourier series. The random Fourier
series also converge or diverge at almost every point with
probability one depending on whether the coefficients are square
summable (see, e.g. \cite{Kahane}). A similar result holds for the
lacunary Fourier series (see, e.g. \cite{Zygmund} for further
references). This analogy is not accidental. Indeed, the spectral
properties are related to the behavior of solutions of the
Schr\"odinger equation. Although the precise link between
solutions and local (in energy) properties of spectral measure is
given by the subordinacy condition discovered by Gilbert and
Pearson \cite{GP}, boundedness of the solutions is typically
associated with the absolutely continuous spectrum. In particular,
it has been shown by Stolz \cite{Stolz} and by Simon \cite{Sim2}
that if all solutions of the equation $-u'' +V(x)u = E u$ are
bounded for each $E$ in a set $S$ of positive Lebesgue measure,
then the absolutely continuous part of the spectral measure gives
positive weight to $S.$ Establishing the boundedness of the
solutions, on the other hand, can be thought of as a nonlinear
analog of proving the convergence of Fourier series, at least for
the potentials in $L^2(\Rm).$ To clarify this idea, it is
convenient to introduce the generalized Pr\"ufer transform, a very
useful tool for studying the solutions in one dimension. We will
very roughly sketch the idea behind Theorems~\ref{sparse1},
\ref{random1} following \cite{KLS}.

Let $u(x,k)$ be a solution of the eigenfunction equation
\begin{equation}\label{ef}
-u'' +V(x)u = k^2 u.
\end{equation}
 The modified Pr\"ufer variables are
introduced by
\begin{equation}\label{pruf}
u'(x,k) = k R(x,k) \cos \theta(x,k); \,\,\,u(x,k) = R(x,k) \sin
\theta(x,k).
\end{equation}
The variables $R(x,k)$ and $\theta(x,k)$ verify
\begin{eqnarray}\label{pruferam}
(\log R(x,k)^2)' = \frac{1}{k}V(x) \sin 2\theta(x,k) \\
\label{pruferan} \theta(x,k)' = k - \frac{1}{k} V(x) (\sin
\theta)^2.
\end{eqnarray}
Fix some point $x_0$ far enough, and set $\theta(x_0)=\theta_0,$
$R(x_0) = R_0.$ From \eqref{pruferan}, we see that
\begin{equation}\label{anggrow} \theta(x) = \theta_0 +
k(x-x_0) - \frac{V(x)}{k} (\sin (k(x-x_0)+\theta_0))^2 + O(V^2) =
\theta_0 + k(x-x_0) + \delta \theta +O(V^2).
\end{equation}
Then
\[ \sin 2\theta (x,k) = \sin (2(\theta_0 + k(x-x_0))) + 2\cos
(2(\theta_0+k(x-x_0)))\delta \theta + O(V^2). \] From
\eqref{anggrow} and the equation \eqref{pruferam} for the
amplitude, we find
\begin{equation}\label{ampgrow}
\frac{d}{dx} (\log (R^2(x)) = t_1 + t_2 +O(V^3),
\end{equation}
where
\[ t_1 = \frac{V(x)}{k} \sin (2(\theta_0 + k(x-x_0))) -
\frac{V(x)}{2k^2} \left( \int\limits_{x_0}^x V(y)\,dy \right) \cos
(2(\theta_0 + k(x-x_0))), \] and
\[ t_2 = \frac{1}{4k^2} \frac{d}{dx} \left[ \int\limits_{x_0}^x
V(y) \cos (2(\theta_0+k(y-x_0)))\,dy \right]^2. \] In both the
random and sparse cases, we obtain the asymptotic behavior of
$R(x)$ by summing up contributions from finite intervals. In the
random case, these intervals correspond to the independent random
parts of the potential, while in the sparse case $R(x)$ remains
unchanged between the neighboring bumps, and we only have to add
the contributions of the bumps. In both cases, for different
reasons, the contributions of the $t_1$ terms can be controlled
and are finite (with probability one in the random case). For
random potentials, one uses the linearity of terms entering $t_1$
in $V$ and the independence of different contributions; the
argument is then similar to the Fourier transform case and gives
convergence as far as $V \in L^2$ by Kolmogorov three series
theorem. We note that for the second term in $t_1$ one actually
has to use a bit more subtle reasoning also taking into account
the oscillations in energy. In the sparse case, one uses the fact
that contributions from different steps are oscillating in $k$
with very different frequency due to the large distance between
$x_n$ and $x_{n+1}.$ Again, the argument is related to the
techniques used to study the lacunary Fourier series. On the other
hand, the sum of $t_2$ terms is finite if $V \in L^2,$ leading to
the boundedness of the solutions and absolutely continuous
spectrum. If $V$ is not $L^2,$ the sum of $t_2$ terms diverges and
can be shown to dominate the other terms due to lack of sign
cancelations.

The question whether the absolutely continuous spectrum is
preserved for general $L^2$ potentials remained open longer. The
initial progress in this direction focused on proving boundedness
of solutions for almost every energy. There are many examples,
starting from the celebrated Wigner and von Neumann \cite{WvN}
construction of an imbedded eigenvalue, which show that spectrum
does not have to be purely absolutely continuous if $V \notin
L^1,$ and imbedded singular spectrum may occur. We will discuss
some of these examples in Section~\ref{strikex}. Thus there can be
exceptional energies with decaying and growing solutions. Again,
one can think of a parallel with Fourier transform, where the
integral $\int_{-N}^N e^{ikx} g(x)\,dx$ may diverge for some
energies if $g \in L^2.$ It was conjectured by Luzin early in the
twentieth century that nevertheless the integral converges for
a.e. $k.$ The question turned out to be difficult, and required an
extremely subtle analysis by Carleson to be solved positively in
1955 in a famous paper \cite{Car}. If $g \in L^p$ with $p<2,$ the
problem is significantly simpler, and has been solved by Zygmund
in 1928 \cite{Zyg} (see also Menshov \cite{Men} and Paley
\cite{Pal} for the discrete case). As the equation
\eqref{pruferam} suggests, the problem of boundedness of solutions
to Schr\"odinger equation may be viewed as a question about a.e.
convergence of a nonlinear Fourier transform. Research in this
direction started from work of Christ, AK, Molchanov and Remling
on power decaying potentials \cite{Kis1,Kis2,Mol,CK1,Rem2}. An
elegant and simple paper by Deift and Killip \cite{DK} used a
completely different idea, the sum rules, to prove the sharp
result, the preservation of the absolutely continuous spectrum for
$L^2$ potentials. The two approaches can be regarded as
complementary: the study of solutions gives more precise
information about the operator and dynamics, but has so far been
unable to handle the borderline case $p=2.$ The sum rules methods
give the sharp result on the nature of the spectrum, but less
information about the nature of the eigenfunctions and dynamical
properties. We will briefly sketch the most current results in
both areas, starting with the solutions approach.

Let $H_V$ be the whole line Schr\"odinger operator. Recall that
the modified wave operators are defined by
\begin{equation}\label{modwave}
\Omega^\pm_m g = L^2-\lim_{t \rightarrow \mp \infty} e^{iH_Vt}
e^{-iW(-i\partial_x,t)} g,
\end{equation}
where the operator $e^{-iW(-i\partial_x,t)}$ acts as a multiplier
on the Fourier transform of $g.$ Let
\begin{equation}\label{WKB}
W (k,t) = k^2 t + \frac{1}{2k} \int_0^{2kt} V(s)\,ds.
\end{equation}
The following theorem has been proved in \cite{CK2}.
\begin{thm}\label{WKBth}
Assume that $V \in L^p,$ $p<2.$ Then for a.e. $k,$ there exist
solutions $u_\pm (x,k)$ of the eigenfunction equation \eqref{ef}
such that
\begin{equation}\label{wkbas}
 u_\pm(x,k) = e^{ikx - \frac{i}{2k} \int_0^x V(y)\,dy} (1+o(1))
 \end{equation}
as $x \rightarrow \pm \infty.$ Moreover, the modified wave
operators \eqref{modwave} exist.

Assume, in addition, that $V(x)|x|^\gamma \in L^p$ for some $p<2$
and $\gamma>0.$ Then the Hausdorff dimension of the set of $k$
where \eqref{wkbas} fails cannot exceed $1 - \gamma p'$ (where
$p'$ is the H\"older conjugate exponent to $p$).
\end{thm}
The asymptotic behavior in \eqref{wkbas} as well as the phase in
\eqref{modwave} coincide with the WKB asymptotic behavior, which
has been known for a long time for potentials satisfying
additional conditions on the derivatives. The main novelty of
\eqref{modwave} is that no such condition is imposed. Note that if
the integral $\int_0^N V(s)\,ds$ converges, the asymptotic
behavior of $u_\pm$ becomes identical to the solutions of the
unperturbed equation, and modified wave operators can be replaced
by the usual M\"oller wave operators. The proof of
Theorem~\ref{WKBth} proceeds by deriving an explicit series
representation for the solutions $u_\pm$ via an iterative
procedure. The terms in the series may diverge for some values of
$k,$ but converge almost everywhere. The first term in the series
is a generalization of the Fourier transform, $\int_0^N \exp(ikx -
\frac{i}{k} \int_0^x V(y)\,dy) V(x)\,dx.$ The main difficulty in
the proof comes from proving the estimates for the multilinear
higher order terms such that the series can be summed up for a.e.
$k.$ The estimate \eqref{wkbas} implies that all solutions of
\eqref{ef} are bounded for a.e. $k$ if $V \in L^p,$ $p<2,$ and can
be thought of as a nonlinear version of Zygmund's result for the
Fourier transform. The question whether \eqref{wkbas} holds and
whether the modified operators exist for $V \in L^2$ is still
open, and appears to be very hard, especially the a.e. boundedness
of the eigenfunctions. Indeed, proving \eqref{wkbas} would be the
nonlinear analog of the Carleson theorem. Moreover, Muscalu, Tao
and Thiele showed \cite{MTT} that the method of \cite{CK2} has no
chance of succeeding when $p=2$ (since some terms in the
multilinear series expansion may diverge on a set of positive
measure). The techniques behind Theorem~\ref{WKBth} have been used
to prove related results on slowly varying potentials (with
derivatives in $L^p,$ $p<2$) and perturbations of Stark operators.
See \cite{CK3,CK4} for more details.

The sum rules approach to proving absolute continuity of the
spectrum has been pioneered by Deift and Killip and led to an
explosion of activity in the area and many impressive new results.
Assume that $V(x) \in C_0^\infty(\mathbb{R}).$ Let us consider the
solution $f(x,k)$ of \eqref{ef} such that $f(x,k) =\exp(ikx)$ as
$x$ is to the right of the support of $V$. Then, for $x$ to the
left of the support of $V$, we have
$f(x,k)=a(k)\exp(ikx)+b(k)\exp(-ikx)$. The solution $f(x,k)$ is
called the Jost solution, and $f(0,k)$ the Jost function.
Coefficient $t(k)=a^{-1}(k)$ is the transmission coefficient in
classical scattering theory. Denote $E_j$ the eigenvalues of the
operator $H_V.$ The following identity is well known (see e.g.
\cite{FT}):
\begin{equation}\label{sr1}
\int\limits_{-\infty}^\infty (\log |a(k)|) k^2\,dk + \frac{2\pi}{3}
\sum\limits_j |E_j|^{3/2} =
\frac{\pi}{8}\int\limits_{-\infty}^\infty V^2(x)\,dx.
\end{equation}
The identity can be proved, for example, by a contour integration in
the complex upper half plane of an asymptotic expansion in $k^{-1}$
of the integral equation one can write for $f(x,k)$. There is a
whole hierarchy of formulas (sum rules) similar to \eqref{sr1}. This
fact is related to the role that the inverse scattering transform
for Schr\"odinger operators plays in understanding the KdV dynamics.
The expressions involving $V(x)$ which appear on the right hand side
in such sum rules are the KdV invariants. The inequalities of type
\eqref{sr1} have been applied in the past to derive bounds on the
moments of the eigenvalues of $H_V$ (Lieb-Thirring inequalities).
Deift and Killip realized that the coefficient $a(k)$ is directly
linked to the spectral measure of $H_V.$ Building a sequence of
compact support approximations to $V(x),$ and then passing to the
limit, one essentially derives a lower bound on the entropy of the
absolutely continuous part of the spectral measure:
\[ \int\limits_{I} \log \mu'(\lambda) d\lambda >
-\infty, \] where $d\mu(\lambda)$ is the spectral measure and $I$
is an arbitrary bounded subinterval of $\mathbb{R}^+$. This proves
\begin{thm}\label{dk}
For any $V \in L^2,$ the essential support of the absolutely
continuous spectrum of the operator $H_V$ coincides with the half
axis $(0,\infty).$ That is, the absolutely continuous part of the
spectral measure $\mu_{ac}$ gives positive weight to any set $S
\subset (0,\infty)$ of the positive Lebesgue measure.
\end{thm}
Killip \cite{Kil} later proved a strengthened version of
Theorem~\ref{dk}, also applicable to potentials from $L^3$ given
additional assumptions on the Fourier transform, and to Stark
operators. The key advance in \cite{Kil} is a local in energy
version of \eqref{sr1}, which is more flexible and useful in
different situations. The important fact exploited in \cite{Kil}
is that the Jost function is actually the perturbation determinant
of the Schr\"odinger operator. That yields a natural path to
obtaining estimates necessary to control the boundary behavior of
Jost function. The square of the inverse of the Jost function, on
the other hand, is proportional to the density of the spectral
measure (see \eqref{factor} for a similar higher dimensional
relation). Therefore, the estimates on Jost function have deep
spectral consequences.

The results of \cite{DK} have been extended to slowly varying
potentials with higher order derivative in $L^2$ by Molchanov,
Novitski and Vainberg \cite{MNV}, using the higher order KdV
invariants. Some improvements were made in \cite{den7} where the
asymptotical methods for ODE were used.

 In the discrete setting, the
application of sum rules led Killip and Simon \cite{KS} to a
beautiful result giving a complete description of the spectral
measures of Jacobi matrices which are Hilbert-Schmidt perturbations
of a free Jacobi matrix. Further extensions to slower decaying
perturbations of Jacobi matrices and Schr\"odinger operators have
been obtained in different works by Laptev, Naboko, Rybkin and
Safronov,
\cite{LNS1,LNS2,Ryb,Saf}. Recently, Killip and Simon
\cite{KSC} proved a continuous version of their Jacobi matrix
theorem, giving a precise description of spectral measures that
can occur for Schr\"odinger operators with $L^2$ potentials. We
will further discuss their result in the following section.

Certain extensions of the sum rules method also have been applied
to higher dimensional problems, and will be discussed in
Section~\ref{hd}.

We complete this section with a somewhat philosophical remark. The
technique of Deift-Killip proof (and its developments) has a
certain air of magic about it. After all, it is based on the
identity - sum rule \eqref{sr1}, a rarity in analysis. Recall the
classical von Neumann-Kuroda theorem, which says that given an
arbitrary self-adjoint operator $A,$ one can find an operator $Y$
with arbitrary small Hilbert-Schmidt norm (or any Schatten class
norm weaker than trace class) such that $A+Y$ has pure point
spectrum. Theorem~\ref{dk} says that the situation is very
different if one restricts perturbations to potentials in the case
of $A = H_0.$ The result is so clear cut that one has to wonder if
there is a general, operator theory type of result which, for a
given $A$ with absolutely continuous spectrum, describes classes
of perturbations which would be less efficient in diagonalizing
it. Such more general understanding could prove useful in other
situations, but currently is completely missing.

\section{The Striking Examples}\label{strikex}

Apart from the general results described in the previous section,
there are fairly explicit descriptions of decaying potentials
leading to quite amazing spectral properties. The examples we
discuss here deal with imbedded singular spectrum. Although all
examples we mention are constructed in one dimension, in most
cases it is not difficult to extend them to an arbitrary dimension
using spherically symmetric potentials. The grandfather of all
such examples is a Wigner-von Neumann example of a potential which
has oscillatory asymptotic behavior at infinity, $V(x) = 8\sin
2x/x + O(x^{-2}),$ and leads to an imbedded eigenvalue at $E=1.$
The imbedded singular spectrum for decaying potentials is the
resonance phenomenon, and requires oscillation in the potential,
similarly to the divergence of Fourier series or integrals. It is
also inherently unstable - for example, for a.e. boundary
condition in the half line case there is no imbedded singular
spectrum. The first examples we are going to discuss are due to
Naboko \cite{Nab} and Simon \cite{Sim1}, who provided different
constructions for potentials leading to a similar phenomenon.
\begin{thm}\label{dpis}
For any positive monotone increasing function $h(x) \rightarrow
\infty,$ there exist potentials satisfying $|V(x)| \leq
\frac{h(x)}{1+|x|}$ such that the half-line operator $H_V$ (with,
say, Dirichlet boundary condition) has dense point spectrum in
$(0,\infty).$

If $|V(x)| \leq \frac{C}{1+|x|},$ the eigenvalues $E_1, \dots,
E_n, \dots$ of $H_V$ lying in $(0,\infty)$ must satisfy $\sum_n
E_n < \infty.$
\end{thm}
The last statement of Theorem~\ref{dpis} has been proved in
\cite{KLS}.

The construction of Naboko used the first order system
representation of the Schr\"odinger equation, and had a
restriction that the square roots of eigenvalues in $(0,\infty)$
had to be rationally independent. Simon's construction can be used
to obtain any dense countable set of eigenvalues in $(0,\infty).$
The idea of the latter construction is, roughly, given a set of
momenta $k_1,\dots, k_n, \dots,$ take
\[ V(x) = W(x) + \sum_n \chi_{(x_n, \infty)}(x) B_n \frac{\sin (2 k_n x +
\beta_n)}{x}. \] Here $x_n,$ $B_n$ and $\beta_n$ have to be chosen
appropriately, and $W(x)$ is a compactly supported potential whose
job is to make sure that the $L^2$ eigenfunctions at $E_n$ satisfy
the right boundary condition at zero. Thus basically, the
potential is a sum of resonant pieces on all frequencies where the
eigenvalues are planned.
To explain the argument better, we will outline the third
construction of such an example, which in our view is technically
the simplest one to implement. We will only sketch the proof; the
details are left to the interested reader.

\begin{proof}[Proof of Theorem~\ref{dpis}]
Without loss of generality, we assume that $h(x)$ does not grow
too fast, say $|h(x)| \leq x^{1/4}.$ Recall the Pr\"ufer variables
$R(x,k),$ $\theta(x,k)$ and equations \eqref{pruferan},
\eqref{pruferam} they satisfy. For $x \leq x_1,$ $x_1$ to be
determined later, let
\begin{equation}\label{Vex}
 V(x) = -\frac{h(x)}{2(1+|x|)} \sin
2\theta(x,k_1).
\end{equation}
Here $k_1^2$ is the first eigenvalue from the list we are trying
to arrange. Note that the seeming conflict between defining $V$ in
terms of $\theta$ and $\theta$ in \eqref{pruferan} in terms of $V$
is resolved by plugging \eqref{Vex} into \eqref{pruferan}, solving
the resulting nonlinear equation for $\theta(x,k_1),$ and defining
$V(x)$ as in \eqref{Vex}. Now if $V$ is defined according to
\eqref{Vex} on all half-axis, one can see using \eqref{pruferam},
\eqref{pruferan} and integration by parts that
\[ \log (R(x,k_1)^2) = -\int_0^x \frac{h(y)}{2k_1(1+|y|)}\,dy +
O(1). \] Because of our assumptions on $h(x),$ $R(x,k_1)$ is going
to be square integrable.

Now we define our potential $V(x)$ by
\begin{equation}\label{Vex1}
 V(x) = -\sum\limits_{j=1}^{\infty} \frac{h(x)}{2^j(1+|x|)} \chi_{(x_j, \infty)}(x) \sin
2\theta(x,k_j).
\end{equation}
Each $x_n > x_{n-1}$ is chosen inductively, so that the following
condition is satisfied: for any $j<n,$
\begin{equation}\label{nonres}
 {\rm sup}_x \left| \int_{x_n}^x \frac{h(y)}{(1+|y|)} \sin
 2\theta(y,k_n) \sin 2\theta(y,k_j) \,dy \right| \leq 1.
\end{equation}
Using \eqref{pruferan} and integration by parts, it is easy to see
that on each step, the condition \eqref{nonres} will be satisfied
for all sufficiently large $x_n.$ A calculation using
\eqref{pruferan}, \eqref{pruferam}, and integration by parts then
shows that $R(x,k_n)$ is square integrable for each $n.$ From
\eqref{Vex1} it also follows that $V(x) \leq h(x)/(1+|x|).$
\end{proof}

The examples with imbedded singular continuous spectrum are
significantly harder to construct. The main difficulty is that
while to establish point spectrum one just needs make sure that
$L^2$ norm of the solution is finite, it is not quite clear what
does one need to control to prove the existence of the singular
continuous component of the spectral measure. At the ICMP Congress
in London \cite{Simrev}, Simon posed a problem of finding a
decaying potential leading to imbedded singular continuous
spectrum. First progress in this direction has been due to Remling
and Kriecherbauer \cite{Rem3,KR}. In particular, they constructed
fairly explicit examples of potentials satisfying $|V(x)| \leq
C(1+|x|)^{-\alpha}$, $\alpha
>1/2,$ such that the Hausdroff dimension of the set of singular
energies where the WKB asymptotic behavior \eqref{wkbas} fails is
equal to $2(1-\alpha).$ This is sharp according to
Theorem~\ref{WKBth} (and earlier work of Remling \cite{Rem4} on
power decaying potentials). The set of singular energies is the
natural candidate to support the singular continuous part of the
measure, but the actual presence of the singular continuous part
of the measure remained open.

The first breakthrough came in a work of SD \cite{Den2} where the
following was proved
\begin{thm}\label{scl2}
There exist potentials $V \in L^2$ such that the operator $H_V$
has imbedded singular continuous spectrum in $(0,\infty).$
\end{thm}
The method was inspired by some ideas in approximation theory (see
the next section) and by inverse spectral theory. The classical
inverse spectral theory results (see e.g. \cite{Lev,Mar}) imply
that one can find potentials corresponding to spectral measures
with imbedded singular continuous component. The standard
procedure, however, does not guarantee a decaying potential. In
the meantime, one can develop an inverse spectral theory type of
construction where one also controls the $L^2$ norm of the
potentials corresponding to certain approximations of the desired
spectral measure, in the limit obtaining the $L^2$ potential. The
key control of the $L^2$ norm appears essentially from the sum
rule used by Deift and Killip. The construction in \cite{Den2} was
employing some estimates for the Krein systems, a more general
system of first order differential equations. The amazing aspect
of the construction has been that there is a great flexibility on
how the singular part of the spectral measure may look. Later,  in
the paper \cite{DenK}, it was proved that the imbedded singular
continuous spectrum can occur for faster decaying potentials,
namely if
\[ \int_x^{\infty} q^2(t)\,dt \le C(1+x)^{-1+D+\epsilon}, \]
then the spectral measure can have a singular continuous component
of exact dimension $D.$

Killip and Simon \cite{KS} realized that the idea of \cite{Den2}
is not tied to the Krein systems. They proved a comprehensive
theorem, providing a complete characterization of the spectral
measures of Jacobi matrices which are Hilbert-Schmidt
perturbations of the free matrix. This theorem should be regarded
as an analog of the celebrated Szeg\H{o} theorem for polynomials
orthogonal on the unit circle \cite{szego, ba4}. Recently, they
also extended their result to the continuous case, where it reads
as follows \cite{KSC}. Denote $d\rho(E)$ the spectral measure of
$H_V,$ set $d \rho_0(E) = \pi^{-1} \chi_{[0,\infty)}(E)
\sqrt{E}dE,$ and define a signed measure $\nu(k)$ on $(1,\infty)$
by
\[ \frac{2}{\pi} \int f(k^2) k d\nu(k) = \int f(E) [d\rho(E)
-d\rho_0(E)].\]
 Given a (signed) Borel measure $\nu,$ define
\[ M_s\nu (k) = {\rm sup}_{0 < L \leq 1} \frac{1}{2L} |\nu|
([k-L,k+L]). \] Denote $d \mu/d \sigma$ the Radon-Nikodym
derivative of $\mu$ with respect to $\sigma.$

\begin{thm}\label{kscth}
A positive measure $\rho$ on $\Rm$ is the spectral measure
associated to a $V \in L^2(\Rm^+)$ if and only if \\
{\rm (i)} ${\rm supp} (d\rho) = [0,\infty) \cup \{ E_j \}_{j=1}^N$
with
$E_1 <E_2< \dots <0$ and $E_j \rightarrow 0$ if $N= \infty.$ \\
{\rm (ii)} \begin{equation}\label{normcon} \int \log \left[ 1+
\left( \frac{M_s\nu(k)}{k}\right)^2 \right] k^2 \,dk < \infty
\end{equation}
{\rm (iii)} \[ \sum\limits_j |E_j|^{3/2} < \infty \]
{\rm (iv)} \[
\int\limits_0^\infty \log \left[ \frac14 \frac{ d \rho}{d \rho_0}
+ \frac12 +\frac14 \frac{d \rho_0}{d\rho} \right] \sqrt{E}\,dE <
\infty. \]
\end{thm}
The theorem shows explicitly that the singular part of the
spectral measure corresponding to an $L^2$ potential can be pretty
much anything on the positive half-axis, as far as a certain
normalization condition \eqref{normcon} is satisfied.

The last example that we would like to mention provides the sharp
rate of decay for which the imbedded singular continuous spectrum
may appear \cite{Kis4}.
\begin{thm}\label{imsc}
For any positive monotone increasing function $h(x) \rightarrow
\infty,$ there exist potentials satisfying $|V(x)| \leq
\frac{h(x)}{1+|x|}$ such that the half-line operator $H_V$ (with,
say, Dirichlet boundary condition) has imbedded singular
continuous spectrum. The potential $V(x)$ can be chosen so that
the M\"oller wave operators exist, but are not asymptotically
complete due to the presence of the singular continuous spectrum.

On the other hand, if $|V(x)| \leq \frac{C}{1+|x|},$ the singular
continuous spectrum of $H_V$ is empty.
\end{thm}
The proof is based on building a sequence of approximation
potentials $V_n$ which have, respectively,  $2^n$ imbedded
eigenvalues $E^n_j$, approaching a Cantor set. The key is to
obtain uniform control of the norms of the corresponding
eigenfunctions, $\|u(x,E^n_j)\|^2_2 \leq C2^{-n}.$ Such estimate
allows to control the weights the spectral measure assigns to each
eigenvalue, and to pass to the limit obtaining a non-trivial
singular continuous component. The estimate of the norms of the
eigenfunctions is difficult and is proved using the Pr\"ufer
transform, and a Splitting lemma allowing to obtain two imbedded
eigenvalues from one. This lemma is based on a model nonlinear
dynamical system providing an elementary block of construction.

\section{Dirac operators, Krein systems, Jacobi matrices, and OPUC}

It was understood a long time ago that the spectral theory of
one-dimensional differential operators (Schr\"odinger, Dirac,
canonical systems) has a lot in common with the classical theory
of polynomials orthogonal on the real line. These polynomials are
eigenfunctions of the Jacobi matrix, also quite classical object
in analysis. So, naturally, to understand the problems for
differential operators one might first study analogous problems
for the discrete version. Unfortunately, the Jacobi matrices are
not so easy to study as well. That difficulty was encountered by
many famous analysts (such as Szeg\H{o}) and the answer was found
in the theory of polynomials orthogonal on the unit circle. It
turns out that for many questions (especially in scattering
theory) that is more natural and basic object to study. Then many
results and ideas can be implemented for Jacobi matrices. For
differential operators, the situation is similar. In many cases,
instead of Schr\"odinger operator it makes sense to consider Dirac
operator and for a good reason. Already in 1955 M. Krein
\cite{krein} gave an outline of the construction that led to the
theory of continuous analogs of polynomials orthogonal on the
circle. Instead of complex polynomials, one has the functions of
exponential type that satisfy the corresponding system of
differential equations (the Krein system)
\begin{equation}
\left\{
\begin{array}{ll}
P^{\prime}(r,\lambda)=i\lambda P(r,\lambda)-\bar A(r) P_*(r,\lambda), & P(0,\lambda)=1, \\
P_*^{\prime}(r,\lambda)=-A(r)P(r,\lambda), & P_*(0,\lambda)=1%
\end{array}
\right.  \label{krein2}
\end{equation}
Although more complicated than the OPUC case, the corresponding
theory can be developed. It turns out that the Krein systems happen
to be in one-to-one correspondence with the canonical Dirac
operators:

\begin{equation} {D}\left[
\begin{array}{c}
f_{1} \\
f_{2}%
\end{array}%
\right] =\left[
\begin{array}{cc}
-b & d/dr-a \\
-d/dr-a & b%
\end{array}%
\right] \left[
\begin{array}{c}
f_{1} \\
f_{2}%
\end{array}%
\right], f_2(0)=0  \label{Dir}
\end{equation}%
In fact, $a(r)=2\Re {A(2r)}, b(r)=2\Im {A(2r)}$.

In his pioneering paper \cite{krein}, Krein states the following
\begin{thm}\label{pde1}
If $a(r), b(r)\in L^2(\mathbb{R^+})$ then
$\sigma_{ac}(D)=\mathbb{R}$. Also,
\[
\int\limits_{-\infty}^\infty \frac{\ln
\sigma'(\lambda)}{1+\lambda^2}>-\infty
\]
where $\sigma$ is the spectral measure for Dirac operator.
\end{thm}
This result is actually the corollary from the analogous statement
for Krein systems. Just like  OPUC are related to Jacobi matrices,
the Krein systems and Dirac operators generate Schr\"odinger
operators. To be more accurate, the study of Schr\"odinger operator
with decaying potential is essentially equivalent to the study of
Dirac operator with $b=0, a=C+q$ where the constant $C$ is large
enough and $q$ decays at infinity in essentially the same way as the
Schr\"odinger potential does \cite{den3, den7}.

 We already mentioned the problem of proving the existence of wave
operators for $L^2$ potentials in the second section. In
\cite{Den1}, it was proved that for the Dirac operator, wave
operators do exist if $a, b \in L^2.$ The proof bypasses the
question about a.e. in energy behavior of the eigenfunctions, and
employs instead integral estimates. The key difference between
Dirac and Schr\"odinger cases is different free evolution. One
manifestation of this difference is the fact that no WKB
correction is needed in the definition of wave operators; usual
M\"oller wave operators exist for the $L^2$ perturbations of the
Dirac operator. Nevertheless, the result may indicate that the
$L^2$ wave operator question for Schr\"odinger operators is easier
to resolve than the question of the asymptotic behavior
\eqref{wkbas}.

The study of Dirac operators is often more streamlined than
Schr\"odinger in both one-dimensional and multidimensional cases
(see the next section), but it already poses significant technical
difficulties whose resolution is far from obvious and proved to be
very fruitful for the subject in general.

\section{The Multidimensional Case}\label{hd}

As opposed to one-dimensional theory, the spectral properties of
Schr\"odinger operators with slowly decaying potentials in higher
dimensions are much less understood. Early efforts focused on the
short range case, $|V(x)| \leq C(1+|x|)^{-1-\epsilon},$ culminating
in the proof by Agmon \cite{Ag} of the existence and asymptotic
completeness of wave operators in this case. In the long range case,
H\"ormander \cite{Hor} considered a class of symbol like potentials,
proving existence and completeness of wave operators. For a review
of these results and other early literature, see \cite{RS3, Yafaev}.
For potentials with less regular derivatives, the conjecture by
Barry Simon \cite{Simrev} states that the absolutely continuous
spectrum of the operator $H_V$ should fill all positive half axis if
\begin{equation}\label{simcon}
\int\limits_{\Rm^d} V^2(x)(1+|x|)^{-d+1}\,dx < \infty.
\end{equation}
To avoid problems with definition of the corresponding self-adjoint
operator (that might appear for dimension high enough because of the
local singularity of potential), we also assume that $V$ belongs to,
say, Kato class: $V\in K_d(\mathbb{R}^d)$ \cite{CFKS}.

Recall that there exist potentials $W(r)$ in one dimension which
satisfy $|W(r)| \leq C r^{-1/2}$ and lead to purely singular
spectrum \cite{KU,KLS}. By taking a spherically symmetric
potential $V(x) = W(|x|),$ one can obtain multidimensional
examples showing that \eqref{simcon} is sharp in many natural
scales of spaces. Notice also that the potential satisfying
(\ref{simcon}) does not have to decay at infinity pointwise in all
directions: it can even grow along some of them. Nevertheless, it
does decay in the average and that makes the conjecture plausible.

Motivated by Simon's conjecture, much of the recent research
focused on long range potentials with either no additional
conditions on the derivatives, or with weaker conditions than in
the classical H\"ormander work. The solutions method so far had
little success in higher dimensions. There are results linking the
behavior of solutions and spectrum which work in higher
dimensions, such as for a example the following theorem proved in
\cite{KL}. In higher dimensional problems, there is no canonical
spectral measure, and the spectral multiplicity can be infinite.
Given any function $\phi \in L^2(\mathbb{R}^d),$ denote $\mu^\phi$
the spectral measure of $H_V$ corresponding to $\phi,$ that is, a
unique finite Borel measure such that $\langle f(H_V)\phi,\phi
\rangle = \int f(E) d\mu^\phi(E)$ for all continuous $f$ with
finite support.
\begin{thm}\label{solutions}
Assume that the potential $V$ is bounded from below. Suppose that
there exists a solution $u(x, E)$ of the generalized eigenfunction
equation
 $(H_{V}-E)u(x,E)=0$
such that
\begin{equation}
\label{1.2} \liminf_{R \to \infty} R^{-1}\int\limits_{|x|\leq R}|
u(x,E)|^{2}\,dx < \infty.
\end{equation}
Fix some $\phi(x) \in C_{0}^{\infty}(\mathbb{R}^d)$ such that
\[ \int\limits_{\mathbb{R}^d}\phi(x) u(x,E)\,dx \neq 0. \]
Then we have
\begin{equation}\label{smest} \limsup_{\delta
\rightarrow 0} \frac{\mu^{\phi}(E-\delta, E+\delta)}{2\delta} >0.
\end{equation}
\end{thm}
Notice that if \eqref{smest} holds on some set $S$ of positive
Lebesgue measure, this implies that the usual Lebesgue derivative
of $\mu^\phi$ is positive a.e. on $S,$ and so presence of the
absolutely continuous spectrum. The condition \eqref{1.2}
corresponds to the power decay $|x|^{(1-d)/2},$ such as spherical
wave solutions decay for the free Laplacian. One may ask how
precise this condition is, perhaps existence of just bounded
solutions on a set $S$ of positive Lebesgue measure is sufficient
for the presence of the absolutely continuous spectrum? It turns
out that, in general, the condition \eqref{1.2} cannot be relaxed.
\begin{theorem}\label{bounsol}
There exists a potential $V$ such that for any $\sigma>0$ there
exists an energy interval $I_\sigma$ with the following
properties.
\begin{itemize}
\item For a.e. $E \in I_\sigma,$ there exists a solution
$u(x,E)$ of the generalized eigenfunction equation satisfying
$|u(x,E)| \leq C(E)(1+|x|)^{\sigma+(1-d)/2}.$
\item The spectrum on $I_\sigma$ is purely singular.
\end{itemize}
\end{theorem}
One way to prove Theorem~\ref{bounsol} is to use one-dimensional
random decaying potentials with $|x|^{-1/2}$ rate of decay. The
results of \cite{KU} or \cite{KLS} show that the spectrum is
singular almost surely, and the eigenfunctions decay at a power
rate. Taking spherically symmetric potentials of this type in
higher dimensions, it is not difficult to see that one gets
examples proving Theorem~\ref{bounsol}.

The link between the behavior of solutions and spectral measures
has been made even more general, sharp and abstract in \cite{CKL}.
However, the difficulty is that obtaining enough information about
solutions in problems of interest is hard: there seems to be no
good PDE analog for the ODE perturbation techniques which can be
used to understand the solutions in one dimension.
On the other hand, the sum rules in higher dimensions typically
involve spaces of potentials which are far from the conjectured
class \eqref{simcon}. Some important progress, however, has been
made using the one-dimensional ideas \cite{Den3,Den4,
Den5,LNS3,LNS4,Perel,Saf1}. After reviewing these results, we will
discuss random decaying potentials in higher dimensions
\cite{B1,B2,RS} as well as quickly mention some interesting recent
progress on a new class of short range potentials \cite{IS} and
imbedded eigenvalues \cite{IJ}.


The OPUC, the Krein systems, and the Dirac operators with
matrix-valued and even operator-valued coefficients were studied
relatively well. This matrix-valued case can give some clues to
the understanding of partial differential equations. Indeed,
writing up, say, Schr\"odinger operator in the spherical
coordinates one obtains the one-dimensional Schr\"odinger operator
with an operator-valued potential. The difficulty is that this
potential $\tilde{V}$ is not bounded since it involves
Laplace-Beltrami operator on the unit sphere. Also
$\tilde{V}(r_1)$ and $\tilde{V}(r_2)$ do not commute for different
values of $r$.

Consider the three dimensional Dirac operator with the following
type of interactions
\begin{equation}
D=-i\alpha \cdot \bigtriangledown +V(x)\beta. \label{dirac8}
\end{equation}
 Here
\begin{equation*}
\alpha _{j}=\left(
\begin{array}{cc}
0 & \sigma _{j} \\
\sigma _{j} & 0%
\end{array}%
\right) ,\ \sigma _{1}=\left(
\begin{array}{cc}
0 & 1 \\
1 & 0%
\end{array}%
\right) ,\ \sigma _{2}=\left(
\begin{array}{cc}
0 & -i \\
i & 0%
\end{array}%
\right) ,\ \sigma _{3}=\left(
\begin{array}{cc}
1 & 0 \\
0 & -1%
\end{array}%
\right)
\end{equation*}%
Matrices $\sigma _{j\text{ }}$ are called the Pauli matrices and
\begin{equation*}
\beta =\left[
\begin{array}{cccc}
1 & 0 & 0 & 0 \\
0 & 1 & 0 & 0 \\
0 & 0 & -1 & 0 \\
0 & 0 & 0 & -1%
\end{array}%
\right]
\end{equation*}
Then the multidimensional analog of Theorem \ref{pde1} says
\cite{Den3} that if $V(x)\in L^\infty$ and the estimate
(\ref{simcon}) holds then $\sigma_{ac}(D)=\mathbb{R}$. Thus the
multidimensional result for Dirac operator is quite satisfactory,
and the Dirac analog of Simon's conjecture \eqref{simcon} holds.
We will sketch the ideas behind this result later. For
Schr\"odinger operator, we first state the following interesting
result by Safronov \cite{Saf1}.

\begin{thm}\label{pde2}
Let $d\geq 3$ and suppose $V\in L^\infty(\mathbb{R}^d)$ is such
that $V(x)\to 0$ as $x\to\infty$. Assume also that $V(x)\in
L^{d+1}(\mathbb{R^d})$ and for some positive $\delta>0$, the
Fourier transform of $V$ satisfies $\hat{V}(\xi)\in
L^2(B_\delta)$. Then, $\sigma_{ac}(H)=\mathbb{R^+}$.
\end{thm}

Note that the $L^{d+1}$ condition corresponds to the
$|x|^{-\frac{d}{d+1}}$ power decay. There are several methods to
obtain that kind of results. The first one, developed by Laptev,
Naboko and Safronov \cite{LNS1}, is based on writing the operator
in the spherical coordinates. Then one takes the Feshbach
projection corresponding to the first harmonic and studies the
corresponding one-dimensional Schr\"odinger operator with
non-local operator-valued potential. Instead of trace equality one
can get an inequality only which still is enough to conclude the
presence of a.c. spectrum. On the other hand, instead of dealing
with Feshbach projections, one can carefully study the
matrix-valued Dirac or Schr\"odinger operator and obtain the
estimates on the entropy of the spectral measure independent of
the size of the matrix \cite{Den3}. Then, an analogous estimate
can be obtained for the corresponding PDE.

Another approach allows to work directly with PDE \cite{Den4}. It
consists of the following observation. Consider, for example,
three-dimensional Schr\"odinger operator with compactly supported
potential. Taking $\phi \in L^\infty$ with compact support we
introduce $u(x,k)=(H-k^2)^{-1}\phi$ with $k\in \mathbb{C}^+, \Re
k>0$. Then, clearly, $u(x,k)$ has the following asymptotic
behavior at infinity:

\[
u(x,k)=\frac{\exp(ikr)}{r}(A_\phi(k,\theta)+\bar{o}(1)), r=|x|,
\theta=x/r
\]
The amplitude $A_\phi(k, \theta)$ can be regarded as an analytic
operator on $L^2(\Sigma)$, where $\Sigma$ is the unit sphere. For
the potential $V$ with compact support, it is continuous up to
each boundary point $k>0$ and the following factorization identity
holds \cite{Yafaev}

\begin{equation}\label{factor}
\mu'_\phi(E)=k\pi^{-1} \|A_\phi(k,\theta)\|^2_{L^2(\Sigma)},\,\,
E=k^2>0.
\end{equation}
Loosely speaking, the density of spectral measure for any vector
$\phi$ can be factorized on the positive interval via the boundary
value of some analytic operator-valued function defined in the
adjacent domain in $\mathbb{C}^+$. Therefore, one can try to
consider the general potential $V$, establish existence of
$A_\phi(k,\theta)$ for all $k\in \mathbb{C}^+, \Re k>0$ with some
(probably rather crude) bounds on the boundary behavior near the
real line. Then the analyticity will be enough to conclude the
necessary estimate on the entropy of the spectral measure,
similarly to the one-dimensional considerations. Here is the
general result.
\begin{thm}\label{pdesol}
Consider potential $V\in L^\infty(\mathbb{R}^d)$. Let
$V_n(x)=V(x)\chi_{|x|<n}$ be its truncation and
$A_{\phi,n}(k,\theta)$-- the corresponding amplitude. Consider an
interval $0<a<k<b$. Assume that for $k_0=(a+b)/2+i\sigma,
\sigma>0$, we have an estimate
$\|A_{\phi,n}(k_0,\theta)\|_{L^2(\Sigma)}>\delta>0$ uniformly in
$n$ and $\|A_{\phi,n}
(\tau+i\varepsilon,\theta)\|_{L^2(\Sigma)}<C\exp(\varepsilon^{-\gamma}),
\gamma>0 $ uniformly in $n, \tau\in [a,b], 0<\varepsilon<2\sigma$.
Then the spectral measure of the function $\phi$ has a.c.
component whose support contains an interval $[a^2, b^2]$.

\end{thm}

The situation is reminiscent of one in the Nevanlinna theory in
the classical analysis when the analyticity and rough bounds close
to the boundary are enough to say a lot about the function. Thus
the whole difficulty here is to obtain the necessary bounds for
the particular PDE. That turns out to be a tricky task but doable
in some cases. For example, the following theorem has been proved
in \cite{Den4}.
\begin{thm}\label{pde4}
Let $Q(x)$ be a $C^1(\mathbb{R}^3)$ vector-field in $\mathbb{R}^3$
and
\[
|Q(x)|<\frac{C}{1+|x|^{0.5+\varepsilon}},\,\,\, |div\;
Q(x)|<\frac{C}{1+|x|^{0.5+\varepsilon}},\,\,\,
|V_1(x)|<\frac{C}{1+|x|^{1+\varepsilon}}, \varepsilon>0,
\]
Then, $H=-\Delta+div\; Q+V_1$ has an a.c. spectrum that fills
$\mathbb{R}^+$.\label{ac-spectrum}
\end{thm}
This theorem is a multidimensional analog of the following result
in dimension one \cite{den3}: if $V=a'+V_1$ where $V_1\in
L^1(\mathbb{R^+})$, $a\in W_2^1(\mathbb{R^+})$, then the a.c.
spectrum of one-dimensional Schr\"odinger with potential $V$
covers the positive half-line. The proof involves
Theorem~\ref{pdesol}, and is based on the uniform estimates for
the Green's function on complex energies. The PDE approach used to
prove Theorem~\ref{pde4} succeeds because it allows not to deal
with negative eigenvalues and the corresponding Lieb-Thirring
inequalities often arising in the sum rule approach (see also
\cite{Kil}).

Although rather elaborate, the conditions on the potential from
the Theorem \ref{pde4} are not very difficult to check.
Essentially, it means that in addition to decay, one has to have
certain oscillation of potential. This condition is related,
although not identical, to the condition on the Fourier transform
in Theorem~\ref{pde2}. For example, an application to random
potentials is possible.
 Consider the following model. Take a smooth function $f(x)$ with the support
inside the unit ball.  Consider

\[
 V_0(x)=\sum\limits_{j\in \mathbb{Z}^+} a_j f(x-x_j)
 \]
 where the points $x_j$ are scattered in
 $\mathbb{R}^3$ such that $|x_k-x_l|>2, k\neq l$, and $a_j\to 0$
 in a way that $|V_0(x)|<C/(1+|x|^{0.5+\varepsilon})$. Let us now
 ``randomize" $V_0$ as follows

 \begin{equation}
 V(x)=\sum\limits_{j\in \mathbb{Z}^+} a_j \xi_j f(x-x_j)
 \label{defin}
 \end{equation}
 where $\xi_j$ are real-valued, bounded, independent random variables with even distribution.
 \begin{thm}\label{pde55}
For $V$ given by (\ref{defin}), we have
$\sigma_{ac}(-\Delta+V)=\mathbb{R^+}$ almost surely.
 \end{thm}
It turns out that for dimension high enough the slowly decaying
random potentials fall into the class considered in Theorem
\ref{pde4}. We note that a similar model has been considered by
Bourgain \cite{B1,B2} in the discrete setting, and will be
discussed below. Theorem~\ref{pde55} and Bourgain's results
suggest very strongly that Simon's conjecture \eqref{simcon} is
true at least in a certain "almost sure" sense (however, the
assumption that the random variables are mean zero is crucial for
the proofs).

The method used in \cite{Den4} was also applied by Perelman in the
following situation \cite{Perel}:

\begin{thm}
For $d=3$,  $\sigma_{ac}(-\Delta+V)=\mathbb{R^+}$ as long as
\[
|V(x)|+|x||\nabla' V(x)|<C/(1+|x|^{0.5+\varepsilon}), \varepsilon>0
\] where $\nabla'$ means the angular component of the gradient.
 \end{thm}
 Here the oscillation is arbitrary in the radial variable and slow
 in the angular variable.
In this case, the Green function has a WKB-type correction which
can be explicitly computed (and has essentially one-dimensional,
integration along a ray, form). The notion of the amplitude
$A_\phi(k,\theta)$ can be modified accordingly and the needed
estimates on the boundary behavior can be obtained.

We now return to the random decaying potentials and discuss the
recent developments in more detail.  Slowly decaying random
potential is a natural problem to tackle if one tries to approach
one of the most important open problems in mathematical quantum
mechanics: the existence of extended states in the Anderson model
in higher dimensions. Important progress in understanding the
random slowly decaying potentials is due to Bourgain \cite{B1,B2}.
Consider random lattice Schr\"odinger operator on $\mathbb{Z}^2$:
 $H_{\omega}=\Delta+V_\omega$, where $\Delta$ is the usual discrete
 Laplacian and $V_\omega$ is a random potential
 \[
 V_\omega=\omega_n v_n
 \]
 with $|v_n|<C|n|^{-\rho}$, $\rho>1/2.$ The random variables $\omega_n$ are
 Bernoulli or normalized Gaussian (and, in particular, are mean zero). Then, \cite{B1}

\begin{thm}\label{bourran} Fix $\tau>0$ and denote $\left. I = \{ E \right|
\tau < |E| < 4-\tau \}.$ Assume that $\rho >1/2,$ and ${\rm sup}_n
|v_n||n^{-\rho}| < \kappa.$ Then for $\kappa = \kappa(\rho,\tau)$
for $\omega$ outside a set of small measure (which tends to zero
as $\kappa \rightarrow 0$) we have \\
1. $H_\omega$ has purely absolutely continuous spectrum on $I$ \\
2. Denoting $E_0(I)$ the spectral projections for $\Delta$ the
wave operators $W_\pm (H,\Delta)E_0(I)$
exists and is complete.
\end{thm}

Using the fact that the absolutely continuous spectrum (and the
existence of the wave operators) are stable under finitely
supported perturbations, one readily obtains absolute continuity
and existence of wave operators almost surely for potentials
satisfying $|v_n| \leq C |n|^{-\rho}.$
The method can also be extended effortlessly to dimensions $d>2.$


Bourgain's approach is based on the careful analysis of the Born's
approximation series for the resolvent. In summation, each of the
terms $[R_0(z)V]^s R_0(z)$ is considered. Then the dyadic
decomposition of $V$ is introduced: $V=\sum_j
V\chi_{2^{j}<|x|<2^{j+1}}$. In the end, the analysis is reduced to
getting the multilinear bounds for the resulting terms. An
interesting and novel in this context ingredient of the proof is
the smart use of a certain entropy bound (the so-called ``dual to
Sudakov" inequality). Later \cite{B2} this approach was further
developed to deal with different situations, such as $L^p$ and
slower power decaying potentials. The main result of \cite{B2} for
the slower power decay is the almost sure existence of a bounded,
not tending to zero solution at a single energy. This, however, is
not yet sufficient for any spectral conclusions. The problem of
handling the random decay with the coefficient $\rho$ even a
little less than $1/2$ remains an interesting open question. So
far, all attempts to deal with this case were not successful.

In another paper on random decaying potentials \cite{RS},
Rodnianski and Schlag showed existence of modified wave operators
for the model with the slow random decay and additional
assumptions ensuring slow variation of the derivatives. The
standard technique of scattering theory, but also with averaging
over the randomness, is employed. That allows to prove scattering
with weaker assumptions than in the standard H\"ormander's case.

Another case for which scattering can be established is the
Schr\"odinger operator on the strip \cite{Den5}. One can show the
presence of a.c. part of the spectrum using the following general
result. Assume that we are given two operators
 $H_1$ and $H_2$ that both act in the same Hilbert space $\mathcal{H}$. Take
$H$ in the following form:
\begin{equation}
H=\begin{bmatrix}
H_1&V\\
V^*&H_2\\
\end{bmatrix} \label{hamil}
\end{equation}
\begin{theorem}\label{n1}
Let $H_1, H_2$ be two bounded self-adjoint operators in the Hilbert
space $\mathcal{H}$. Assume that $\sigma_{ess}(H_2) \subseteq [b, +
\infty]$ and $[a,b] \subseteq \sigma_{ac}(H_1), (a <b)$.  Then, for
any Hilbert-Schmidt $V$, ($V \in {\mathcal{J}}^2$),  we have that
$[a,b] \subseteq \sigma_{ac}(H)$, with $H$ given by (\ref{hamil}).
\end{theorem}

This theorem can be effectively applied to study Schr\"odinger on
the strip. Indeed,
 let
$$
L=-\Delta +Q(x,y),
$$
considered on the strip $\Pi=\{x>0, 0<y<\pi\},$ and impose
Dirichlet conditions on the boundary of $\Pi$.  Consider the
matrix representation of $L$. For $f(x,y)\in L^2(\Pi)$,

$$
f(x,y)=\sqrt{\frac{2}{\pi}}\mathop\sum\limits^{\infty}_{n=1}
\sin(ny) f_n (x),
f_n(x)=\sqrt{\frac{2}{\pi}}\mathop\int\limits_{0}^{\pi} f(x,y)
\sin(ny)dy
$$
and $L$ can be written as follows
\begin{equation}
L=\left[ \begin{array}{ccc}
\displaystyle -\frac{d^2}{dx^2} +Q_{11}(x)+1 & Q_{12}(x) & \ldots\\
Q_{21}(x) & \displaystyle -\frac{d^2}{dx^2}+Q_{22}(x)+4 & \ldots \\
\ldots & \ldots & \ldots \\
\end{array} \right] \label{matr}
\end{equation}
$$
Q_{lj}(x)=\frac{2}{\pi}\mathop\int\limits_0^{\pi}
Q(x,y)\sin(ly)\sin(jy)dy
$$
Assume $\mathop{sup}\limits_{0\leq y \leq \pi}|Q(x,y)|\in L^2 (\Bbb
R^+) \cap  L^\infty (\Bbb{R}^+)$. Since $Q_{11}(x)\in L^2(\Bbb
R^+)$, we can use \cite{DK, MV} and
Theorem \ref{n1} to show the presence of a.c. spectrum.

A great example of how the one-dimensional technique works for
multidimensional problem is the case of scattering on the Bethe
lattice \cite{Den5}. The step-by-step sum rules used by Simon to
study Jacobi matrices \cite{ba} can be adjusted to that case. Let us
consider this model. Take the Caley tree (Bethe lattice)
$\mathbb{B}$ and the discrete Laplacian on it
$$
(H_0 u)_n =\mathop\sum\limits_{|i-n|=1} u_i
$$
 Assume, for simplicity, that the degree at
each point (the number of neighbors) is equal to $3$.  It is well
known that
$\sigma(H_0)=[-2\sqrt{2}, 2 \sqrt{2}]$ and the spectrum is purely
absolutely continuous. Let $H=H_0+V$, where $V$ is a potential.
Consider any vertex $O$. It is connected to its neighbors by three
edges. Delete one edge together with the corresponding part of the
tree stemming from it. What is left will be called
$\mathbb{B}_{O}$. The degree of $O$ within $\mathbb{B}_{O}$ is
equal to $2$.   The solution to Simon's conjecture in this case is
given by the following theorem. We denote by the symbol
$m(\mathbb{B})$ the functional space of sequences decaying at
infinity on $\mathbb{B}$.

\begin{theorem}\label{n2}
If $V\in \ell^\infty(\mathbb{B})\cap m(\mathbb{B}_O)$ and
\begin{equation}
\mathop\sum\limits_{n=0}^{\infty} \frac{1}{2^{n}}
\mathop\sum\limits_{x\in \mathbb{B}_{O}, |x-O|=n} V^2(x) < +\infty
\label{ell2}
\end{equation}
then

$$
[-2\sqrt{2}, 2 \sqrt{2}]=\sigma_{ac}(H_{|\mathbb{B}_{O}})\subseteq
\sigma_{ac}(H)
$$
\end{theorem}
The idea of the proof of Theorem~\ref{n2} is based on \cite{ba}.
This result is sharp in the following sense. Take $V(x):
|V(x)|<C|x-O|^{-\gamma}$. Then for any $\gamma>0.5$ the condition
of the theorem is satisfied (just like in $\mathbb{R}^d$). In the
meantime, one can find the spherically symmetric $V$ with slower
decay $0<\gamma<0.5$, such that there will be no absolutely
continuous spectrum at all.

An important role in the proof is played by the following well-known
recursive relation for
$\langle(H_{|{\mathbb{B}_{O}}}-z)^{-1}\delta_O,\delta_O\rangle$,
where $\delta_O$ is the discrete delta-function at the point $O.$ In
particular, one can derive the following important and physically
meaningful identity

\begin{equation}
\frac{1}{\pi} \mathop\int\limits_{-2\sqrt 2}^{2\sqrt 2}
\sqrt{8-\lambda^{2}} \ln [ \mu'_O(\lambda)] d\lambda \geq
\frac{1}{\pi} \mathop\int\limits_{-2\sqrt 2}^{2\sqrt2}
\sqrt{8-\lambda^{2}} \ln \Big( \frac{\mu'_{O_1}(\lambda) +
\mu'_{O_2}(\lambda)}{2}\Big) d\lambda-V^2(O) \label{thepoint}
\end{equation}
where $\mu'_{O_{1(2)}}$ correspond to the densities of the spectral
measures at points $O_{1(2)}$ on the trees obtained from
$\mathbb{B}_O$ by throwing away the point $O$ along with the
corresponding two edges.  Using inequality $(x+y)/2\geq \sqrt{xy}$
and iterating (\ref{thepoint}), one proves Theorem~\ref{n2}. Formula
(\ref{thepoint}) says, in particular, that no matter what happens
along one branch of the tree, the scattering is possible through the
other branch. It is also clear that the presence of some ``bad"
points in the tree (say, points where we have no control over
potential) should not destroy the scattering as long as these points
are rather ``sparse". What is an accurate measurement of this
sparseness? We suggest the following improvement of
Theorem~\ref{n2}. Consider the tree $\mathbb{B}_O$ with potential
$V$ having finite support, i.e. $V(x)=0$ for $|x-O|>R$. Consider all
paths running from the origin $O$ to infinity without
self-intersections. Using diadic decomposition of the real numbers
on the interval $[0,1]$, we can assign to each path the real number
in the natural way. That is one way of coding the points at
infinity. In principle, this map ${F}$ is not bijection, e.g.
sequences $(1,0,0,\ldots)$ and $(0,1,1,\ldots)$ represent the same
real number $0.5$ but different paths. Fortunately, these numbers
have Lebesgue measure zero and will be of no importance for us. Let
us define the following functions
\[
\phi(t)=\sum\limits_{n=1}^\infty V^2(x_n)
\]
where $x_n$ are all vertices of the path representing the point
$t\in [0,1]$. Since the support of $V$ is within the ball of radius
$R$, function $\phi(t)$ is constant on diadic intervals $[j2^{-R},
(j+1)2^{-R}), j=0,1,\ldots, 2^{R}-1$. Notice that ${F}$ not being
bijection cause no trouble in defining $\phi(t)$.

Define the probability measure $dw(\lambda)=(4\pi)^{-1}
(8-\lambda^2)^{1/2}$ on $[-2\sqrt 2, 2\sqrt 2]$, and
$\rho_O=\sigma'_O(\lambda)[\tilde{\sigma}'(\lambda)]^{-1}$, a
relative density of the spectral measure at the point $O$,
\mbox{($\tilde{\sigma}'(\lambda)=4^{-1}(8-\lambda^2)^{1/2}$).}
Define
\[
s_O=\mathop\int\limits_{-2\sqrt 2}^{2\sqrt 2}  \ln
 \rho_O(\lambda) dw(\lambda)
\]
Consider the probability space obtained by assigning to each path
the same weight (i.e. as we go from $O$ to infinity, we toss the
coin at any vertex and move to one of the neighbors further from $O$
depending on the result). Our goal is to prove the following
\begin{theorem}
For any $V$ bounded on $\mathbb{B}_O$, the following inequality is
true
\begin{equation}
\exp s_O\geq  \mathbb{E} \left\{ \exp \left[-\frac 14
\sum\limits_{n=1}^\infty V^2(x_n)\right]\right\}=\int\limits_0^1
\exp \left(-\frac 14 \phi(t)\right) dt \label{pearson}
\end{equation}
where the expectation is taken with respect to all paths $\{x_n\}$
going from $O$ to infinity without self-intersections. In
particular, if the r.h.s. of (\ref{pearson}) is positive, then
$[-2\sqrt 2,2\sqrt 2]\subseteq \sigma_{ac}(H)$. \label{improve}
\end{theorem}
\begin{proof}

Assume that $V$ has finite support.
 The estimate (\ref{thepoint}) can be rewritten as

\begin{equation}
 \mathop\int\limits_{-2\sqrt 2}^{2\sqrt 2}  \ln
 \rho_O(\lambda) dw(\lambda) \geq
\mathop\int\limits_{-2\sqrt 2}^{2\sqrt2} \ln \Big(
\frac{\rho_{O_1}(\lambda) + \rho_{O_2}(\lambda)}{2}\Big)
dw(\lambda)-V^2(O)/4 \label{thepoint1}
\end{equation}
Now, we will use the Young's inequality
\[
\frac{x^p}{p}+\frac{y^q}{q}\geq xy;\; x, y\geq 0, 1\leq p\leq
\infty, p^{-1}+q^{-1}=1
\]
in (\ref{thepoint1}) to obtain
\begin{equation}
s_O\geq p^{-1}s_{O_1}+q^{-1}s_{O_2}+p^{-1}\ln p+q^{-1}\ln q-\ln
2-V^2(O)/4
\end{equation}
Considering $s_{O_1}$ and $s_{O_2}$ to be fixed parameters and
maximizing the r.h.s. over $p\in [1,\infty]$, we get the following
inequality
\begin{equation}
s_O\geq \ln \frac{\exp s_{O_1}+\exp s_{O_2}}{2}-V^2(O)/4
\label{thepoint2}
\end{equation}
with optimal $p^*=1+\exp(s_{O_2}-s_{O_1})$. Iterate
(\ref{thepoint2}) till we leave the support of $V$. Thus, we get
(\ref{pearson}). Consider now the general case of bounded $V$.
Define the truncation of $V$ to the ball of radius $n$:
$V_n(x)=V(x)\chi_{\{|x-O|\leq n\}}$. For the corresponding
$s^{(n)}_{O}$, we use (\ref{pearson}), take $n\to \infty$ and apply
to the l.h.s. a standard by now argument on the semicontinuity of
the entropy (\cite{KS}, p. 293). Notice that the functions
$\phi^{(n)}(t)$ are nonnegative and increasing in $n$ for each $t$
(this is why $\phi(t)$ is always well defined). Therefore, we get
(\ref{pearson}) from the theorem on monotone convergence.
\end{proof}

Using the Jensen's inequality, we obtain
\begin{corollary}
Assume that $A$ is any subset of $[0,1]$ of positive Lebesgue
measure and $\phi(t)\in L^1(A)$, then
\[
s_O \geq -\frac{1}{4|A|} \int\limits_A \phi(t)dt +\ln|A|
\]
In particular, $[-2\sqrt 2,2\sqrt 2]\subseteq \sigma_{ac}(H)$.
\end{corollary}
It is interesting that the set $A$ does not have to have some
special topological structure, say, to be an interval like in the
standard scattering theory \cite{amrein}.


We now turn to the results on the (generalized) short range
potentials and imbedded eigenvalues in $\mathbb{R}^d$. Recently,
new short-range type results have been obtained for the
multidimensional Schr\"odinger operator with potential from the
$L^p$ and more general classes \cite{GS,IS}. The main goal was to
establish the limiting absorption estimates for the resolvent
acting in certain Banach spaces, which are more detailed and
precise than Agmon's classical results. The standard techniques
developed in the works of Agmon can be improved if one uses the
Stein-Tomas restriction theorem. Here is one of the results in
that direction \cite{IS}:

\begin{thm}
Assume that $V$ is such that
\[
M_{q_0}(V)(x)\in L^{(d+1)/2}(\mathbb{R}^d)
\]
where
\[
M_{q_0}f=\left[\int\limits_{|y|<1/2} |f(x+y)|^q dy\right]^{1/q}
\]
and $q_0=d/2$ if $d\geq 3$, $q_0>1$ if $d=2$. Then, the following is
true for the operator $H=-\Delta+V$:
\begin{itemize}
\item{ The set of nonzero eigenvalues is discrete with the only
possible  accumulation point to be zero. Each nonzero eigenvalue has
finite multiplicity.}
\item{ $\sigma_{s.c.} (H)=\emptyset$, and $\sigma_{a.c.}
(H)=\mathbb{R}^+$}
\item{The wave operators $\Omega^{\pm} (H,H_0)$ exist and are
complete.}
\end{itemize}
\end{thm}
The actual result is a bit stronger. The authors of \cite{IS}
present the whole class of ``admissible" perturbations for which
their method works, including some first order differential
operators.

Another direction in which there has been significant recent
progress concerns imbedded eigenvalues.
The following result, in particular, has been proved by Ionescu
and Jerison \cite{IJ}.
\begin{thm}
Let $V(x)\in L^{d/2}(\mathbb{R}^d)$ and for $d=2$ we also assume
$V(x)\in L^r_{loc}(\mathbb{R}^2), r>1$. Then, $H=-\Delta+V$ does not
have positive eigenvalues. \label{i-j}
\end{thm}
The paper \cite{IJ} actually contains a more general result, which
allows for slower decay of the potential if its singularities are
weaker. The method relied on the Carleman inequality of special
type.

Surprisingly, \cite{IJ} also provides an example of potential $V$
satisfying
\begin{equation}\label{WvNmult}
 |V(x)|<C(|x_1|+x_2^2+\ldots+x_d^2)^{-1}
\end{equation}
for which the positive eigenvalue appears (multidimensional analog
of Wigner- von Neumann potential). From the point of view of
physical intuition, the existence of imbedded eigenvalue for such
potential may seem strange. Indeed, one would expect that tunneling
in the directions $x_2,\dots,x_d$ of fast decay should make the
bound state impossible. Yet, Wigner-von Neumann type oscillation and
slow Coulomb decay in just one direction turn out to be sufficient.
The corresponding eigenfunction decays rather slowly but enough to
be from $L^2$. We note that the potential satisfying \eqref{WvNmult}
just misses $L^{(d+1)/2}(\mathbb{R}^d).$ Quite recently, Koch and
Tataru \cite{KT} improved Theorem \ref{i-j} and showed the absence
of embedded eigenvalued for the optimal $L^{(d+1)/2}(\mathbb{R}^d)$
case. They also considered various long-range potentials and more
general elliptic operators.

There remain many interesting and important open problems
regarding the multidimensional slowly decaying potentials. Simon's
conjecture \eqref{simcon} remains open, and new ideas are clearly
needed to make progress. Improving understanding of random slowly
decaying potentials is another quite challenging direction. Other
difficult and intriguing open questions involve multidimensional
sparse potentials, appearance of imbedded singular continuous
spectrum, and decaying potentials with additional structural
assumptions. This vital area is bound to challenge and inspire
mathematicians for years to come.

{\bf Acknowledgement} \rm The work of AK has been supported in part
by the Alfred P. Sloan Research Fellowship and NSF-DMS grant
0314129. SD was supported by the NSF-DMS grant 0500177.


\begin{thebibliography}{99}



\bibitem{Ag} S.~Agmon, \it Spectral properties of Schr\"odinger operators
and scattering theory, \rm Ann Scuola Norm. Sup. Pisa II, {\bf
2}(1975), 151--218

\bibitem{amrein} W.O.~Amrein, D.B.~Pearson, \it Flux and scattering
 into cones for long range and singular potentials,\rm J. Phys. A, {\bf 30} (1997), no. 15, 5361--5379

\bibitem{BA} M.~Ben-Artzi, \it On the absolute continuity of Schr\"odinger
operators with spherically symmetric, long-range potentials. I,
II. \rm  J. Differential Equations {\bf 38} (1980), 41--50, 51--60

\bibitem{B1} J.~Bourgain, \it On random Schr\"odinger operators on $\Bbb Z\sp 2,$
\rm Discrete Contin. Dyn. Syst. {\bf 8} (2002), 1--15

\bibitem{B2} J.~Bourgain, \it Random lattice Schr\"odinger operators with decaying potential:
some higher dimensional phenomena, \rm Geometric aspects of
functional analysis, Lecture Notes in Math., {\bf 1807} (2003),
Springer, Berlin, 70--98

\bibitem{Car} L.~Carleson,
{\em On convergence and growth of partial sums of Fourier series},
Acta Math.\ {\bf 116} (1966), 135-157

\bibitem{CK1} M.~Christ and A.~Kiselev, \it Absolutely continuous spectrum for one-dimensional
Schr\"odinger operators with slowly decaying potentials: some
optimal results, \rm J. Amer. Math. Soc. {\bf 11} (1998), 771--797

\bibitem{CK2} M.~Christ and A.~Kiselev, \it Scattering and wave operators for one-dimensional
Schr\"odinger operators with slowly decaying nonsmooth potentials,
\rm Geom. Funct. Anal. {\bf 12} (2002), 1174--1234

\bibitem{CK3} M.~Christ and A.~Kiselev, \it WKB and spectral analysis
of one-dimensional Schr\"odinger operators with slowly varying
potentials, \rm Comm. Math. Phys. {\bf 218} (2001), 245--262

\bibitem{CK4} M.~Christ and A.~Kiselev, \it Absolutely continuous
spectrum of Stark operators, \rm Ark. Mat. {\bf 41} (2003), 1--33

\bibitem{CKL} M.~Christ, A.~Kiselev and Y.~Last, {\it
Approximate eigenvectors and spectral theory,} in Differential
Equations and Mathematical Physics, R.~Weikard and G~.Weinstein,
Eds, American Mathematical Society, International Press 2000

\bibitem{CFKS} H.L.~Cycon, R.G.~Froese, W.~Kirsch, B.~Simon, \it Schr\"odinger
operators with application to quantum mechanics and global geometry.
\rm  Texts and Monographs in Physics. Springer Study Edition.
Springer-Verlag, Berlin, 1987.


\bibitem{DSS} F.~Delyon, B.~Simon and B.~Souillard, \it From pure point
to continuous spectrum in disordered systems, \rm Ann. Inst. Henri
Poincare {\bf 42} (1985), 283--309

\bibitem{DK} P.~Deift and R.~Killip, \it On the absolutely continuous spectrum of
one-dimensional Schr\"odinger operators with square summable
potentials, \rm Comm. Math. Phys. {\bf 203} (1999), 341--347

\bibitem{Den1} S.~Denisov, \it On the existence of wave operators for
some Dirac operators with square summable potentials, \rm Geom.
Funct. Anal. {\bf 14} (2004), 529--534

\bibitem{Den2} S.~Denisov, \it On the coexistence of absolutely continuous
and singular continuous components of the spectral measure for
some Sturm-Liouville operators with square summable potential, \rm
J. Differential Equations {\bf 191} (2003), 90--104

\bibitem{Den3} S.~Denisov, \it On the absolutely continuous spectrum of Dirac operator,
\rm Comm. Partial Diff. Eq. {\bf 29} (2004), 1403--1428

\bibitem{Den4} S.~Denisov, \it Absolutely continuous spectrum of multidimensional
Schr\"odinger operator, \rm Int. Math. Res. Not. {\bf 74} (2004),
3963--3982

\bibitem{Den5} S.~Denisov, \it On the preservation of absolutely
continuous spectrum for Schr\"odinger operators, to appear in J.
Funct. Anal.\rm

\bibitem{den3} S.~Denisov, \it On the application of some of M. G. Krein's
results to the spectral analysis of Sturm-Liouville operators, \rm
J. Math. Anal. Appl.  {\bf 261}  (2001),  177--191.

\bibitem{den7} S.~Denisov, \it On the existence of the absolutely continuous component for the
measure associated with some orthogonal systems, \rm Comm. Math.
Phys., {\bf 226}, (2002), 205--220.

\bibitem{DenK} S.~Denisov and S.~Kupin, \it On the singular spectrum
of Schr\"odinger operators with decaying potential, \rm Trans. Amer.
Math. Soc., {\bf 357}, (2005), 1525--1544

\bibitem{D1} A.~Devinatz,  \it The existence of wave operators for oscillating potentials,
\rm  J. Math. Phys. {\bf 21}, (1980), 2406--2411

\bibitem{FT} L.~Faddeev and L.~Takhtajan, \it Hamiltonian methods in the
theory of solitons, \rm Translated from the Russian by A. G. Re\v
iman, Berlin: Springer, 1987

\bibitem{GP} D.J.~Gilbert and D.B.~Pearson, \it On subordinacy and
 analysis
of the spectrum of one-dimensional Schr\"odinger operators, \rm J.
Math. Anal. Appl. {\bf 128} (1987), 30--56

\bibitem{GS} M.~Goldberg and  W.~Schlag, \it A limiting absorption principle for the
three-dimensional Schr\"odinger equation with $L\sp p$ potentials,
\rm Int. Math. Res. Not. (2004), no. 75, 4049--4071

\bibitem{Graf} L.~Grafakos, \it Classical and Modern Fourier
Analysis, \rm Pearson Education, Prentice Hall, 2004

\bibitem{Hor} L.~H\"ormander, \it The existence of wave operators in
scattering theory, \rm Math. Z. {\bf 146} (1976), 69--91

\bibitem{IJ} A.~Ionescu and D.~Jerison, \it On the absence of positive
eigenvalues of Schr\"odinger operators with rough potentials, \rm
Geom. Funct. Anal. {\bf 13} (2003), 1029--1081

\bibitem{IS} A.~Ionescu and W.~Schlag, \it Agmon-Kato-Kuroda theorems
for a large class of perturbations, \rm to appear at Duke
Mathematical Journal

\bibitem{Kahane} J.P.~Kahane, \it Some Random Series of Functions,
\rm Cambridge University Press, Cambridge, 1985

\bibitem{Kil} R.~Killip, \it Perturbations of one-dimensional Schr\"odinger
operators preserving the absolutely continuous spectrum, \rm Int.
Math. Res. Not. {\bf 2002}, 2029--2061

\bibitem{KS} R.~Killip and B.~Simon, \it Sum rules for Jacobi matrices and
their applications to spectral theory, \rm Ann. of Math. {\bf 158}
(2003), 253--321

\bibitem{KSC} R.~Killip and B.~Simon, \it Sum rules and spectral
measures of Schr\"odinger operatros with $L^2$ potentials, \rm
preprint

\bibitem{Kis1} A.~Kiselev \it Absolutely continuous spectrum of one-dimensional Schr\"odinger
operators and Jacobi matrices with slowly decreasing potentials,
\rm Comm. Math. Phys. {\bf 179} (1996), 377--400

\bibitem{Kis2} A.~Kiselev, \it Stability of the absolutely continuous
spectrum of the Schrödinger equation under slowly decaying
perturbations and a.e. convergence of integral operators, \rm Duke
Math. J. {\bf 94} (1998), 619--646

\bibitem{Kis4} A.~Kiselev \it Imbedded singular continuous spectrum for Schr\"odinger
operators, \rm J. Amer. Math. Soc. {\bf 18} (2005), 571--603

\bibitem{KL} A.~Kiselev and Y.~Last, \it Solutions and spectrum
of Schr\"odinger operators on infinite domains, \rm
 Duke Math. J. {\bf 102}(2000), 125--150

\bibitem{KLS} A.~Kiselev, Y.~Last and B.~Simon, \it Modified Pr\"ufer
and EFGP transforms and the spectral analysis of one-dimensional
Schr\"odinger operators, \rm Commun. Math. Phys. {\bf 194}(1998),
1--45

\bibitem{KT} H.~Koch, D.~Tataru, \it Carleman estimates and
absense of embedded eigenvalues, preprint.\rm

\bibitem{KU} S.~Kotani and N.~Ushiroya, \it One-dimensional Schr\"odinger
operators with random decaying potentials, \rm Commun. Math. Phys.
{\bf 115} (1988), 247--266

\bibitem{krein} M.~Krein, \it Continuous analogues of propositions
on polynomials orthogonal on the unit circle, \rm Dokl. Akad. Nauk
SSSR, {\bf 105}  (1955), 637--640.


\bibitem{KR} T.~Kriecherbauer and C.~Remling, \it Finite gap potentials
and WKB asymptotics for one-dimensional Schr\"odinger operators, \rm
Comm. Math. Phys. {\bf 223} (2001), 409--435

\bibitem{LNS1} A.~Laptev, S.~Naboko and O.~Safronov, \it On new relations between
spectral properties of Jacobi matrices and their coefficients, \rm
Comm. Math. Phys. {\bf 241} (2003), 91--110

\bibitem{LNS2} A.~Laptev, S.~Naboko and O.~Safronov, \it Absolutely continuous
spectrum of Jacobi matrices, \rm Mathematical results in quantum
mechanics (Taxco, 2001), 215--223, Contemp. Math., 307, Amer. Math.
Soc., Providence, RI, 2002

\bibitem{LNS3} A.~Laptev, S.~Naboko and O.~Safronov, \it Absolutely continuous spectrum
of Schr\"odinger operators with slowly decaying and oscillating
potentials, \rm Comm. Math. Phys. {\bf 253} (2005), 611--631

\bibitem{LNS4} A.~Laptev, S.~Naboko and O.~Safronov, \it A Szeg\"o condition for a
multidimensional Schr\"odinger operator, \rm J. Funct. Anal. {\bf
219} (2005), 285--305

\bibitem{Lev} B.~Levitan, \it Inverse Sturm-Liouville Problems, \rm
VNU Science Press, Utrecht 1987

\bibitem{Mar} V.~Marchenko, \it Sturm-Liouville Operators and Applications,
\rm Birkhauser, Basel 1986

\bibitem{Men} D.~Menshov, {\em Sur les series de fonctions orthogonales}, Fund.
Math. {\bf 10}, (1927) 375-420

\bibitem{Mol} S.~Molchanov, unpublished

\bibitem{MNV} S.~Molchanov, M.~Novitskii and B.~Vainberg, \it First KdV integrals
and absolutely continuous spectrum for 1-D Schr\"odinger operator,
\rm Comm. Math. Phys. {\bf 216} (2001), 195--213

\bibitem{MV} S.~Molchanov, B.~Vainberg, \it Schr\"odinger operator
with matrix potentials. Transition form the absolutely continuous to
the singular spectrum, \rm J. Funct. Anal. {\bf 215} (2004),
111--129.

\bibitem{MTT} C.~Muscalu, T.~Tao, C.~Thiele, \it A counterexample to a multilinear
endpoint question of Christ and Kiselev, \rm Math. Res. Lett. {\bf
10} (2003), 237--246

\bibitem{Nab} S.N.~Naboko, \it Dense point spectra of Schr\"odinger and
Dirac operators, \rm Theor.-math. {\bf 68} (1986), 18--28

\bibitem{Pal}
R.E.A.C.~Paley, {\em Some theorems on orthonormal functions},
Studia Math. {\bf 3} (1931) 226-245

\bibitem{Pear} D.~Pearson, {\em Singular continuous measures in scattering theory},
Comm. Math. Phys. {\bf 60} (1978), 13--36

\bibitem{Perel} G.~Perelman, \it On the absolutely continous spectrum
of multi-dimensional Schr\"odinger operators, \rm to appear in Int.
Math. Res. Not.

\bibitem{RS3} M.~Reed and B.~Simon, {\it Methods
of Modern Mathematical Physics, III. Scattering Theory}, Academic
Press, London-San Diego, 1979

\bibitem{Rem1} C.~Remling, \it A probabilistic approach to
one-dimensional Schr\"odinger operators with sparse potentials,
\rm Commun. Math. Phys. {\bf 185}(1997), 313--323

\bibitem{Rem2} C.~Remling, \it The absolutely continuous
spectrum of one-dimensional Schr\"odinger operators with decaying
potentials, \rm Comm. Math. Phys. {\bf 193} (1998), 151--170

\bibitem{Rem3} C.~Remling, \it Schr\"odinger operators with
decaying potentials: some counterexamples, \rm Duke Math. J. {\bf
105} (2000), 463--496

\bibitem{Rem4} C.~Remling, \it Bounds on embedded singular
spectrum for one-dimensional Schr\"odinger operators, \rm Proc.
Amer. Math. Soc. {\bf 128} (2000), 161--171

\bibitem{RS} I.~Rodnianski and W.~Schlag, \it
Classical and quantum scattering for a class of long range random
potentials, \rm Int. Math. Res. Not. {\bf 5} (2003), 243--300

\bibitem{Ryb} A.~Rybkin, \it On the spectral $L^2$ conjecture,
$3/2$-Lieb-Thirring inequality and distributional potentials, \rm
preprint

\bibitem{Saf} O.~Safronov, \it The spectral measure of a Jacobi matrix in terms
of the Fourier transform of the perturbation, \rm Ark. Mat. {\bf
42} (2004), 363--377

\bibitem{Saf1} O.~Safronov, \it On the absolutely continuous spectrum of
multi-dimensional Schr\"odinger operators with slowly decaying
potentials, \rm Comm. Math. Phys. {\bf 254} (2005), 361--366

\bibitem{Sim1} B.~Simon, \it Some Schr\"odinger operators with dense point
spectrum, \rm Proc. Amer. Math. Soc. {\bf 125}(1997), 203--208

\bibitem{Sim2} B.~Simon, \it Bounded eigenfunctions and absolutely
continuous spectra for one-dimensional Schr\"odinger operators,
\rm Proc. Amer. Math. Soc. Vol.{\bf 124}(1996), 3361--3369

\bibitem{Simrev} B.~Simon, \it Schr\"odinger operators in the twenty-first century,
\rm Mathematical Physics 2000 (eds. A. Fokas, A. Grigoryan, T.
Kibble and B. Zegarlinski), Imperial College Press, London,
283-288

\bibitem{ba4} B.~Simon, {\it Orthogonal polynomials on the unit
circle}, AMS Colloquium Publications Series, Vol. 54, 2005.


\bibitem{ba} B.~Simon,
\it A canonical factorization for meromorphic Herglotz functions on
the unit disk and sum rules for Jacobi matrices, \rm J. Funct.
Anal., {\bf 214} (2004), 396--409.



\bibitem{Stolz} G.~Stolz, \it Bounded solutions and absolute continuity of
Sturm-Liouville operators, \rm J. Math. Anal. Appl. {\bf 169}
(1992), 210--228

\bibitem{szego}  G.~Szeg\H o, {\em Orthogonal polynomials,} AMS,  Providence,
(1975).


\bibitem{Weid} J.~Weidmann, \it Zur Spektral theorie von Sturm-Liouville
Operatoren, \rm Math. Z. {\bf 98} (1967), 268--302

\bibitem{WvN} J.~von Neumann and E.P.~Wigner, \it \"Uber merkw\"urdige
diskrete eigenwerte - \rm Z. Phys. {\bf 30}(1929), 465--467

\bibitem{White} D.A.W.~White, \it Schr\"odinger operators with rapidly oscillating central potentials,
\rm Trans. Amer. Math. Soc. {\bf 275} (1983), 641--677

\bibitem{Yafaev}
D.~Yafaev, {\em Scattering theory: Some old and new results},
Lecture Notes in Mathematics 1735, Springer, (2000)

\bibitem{Zyg}
A.~Zygmund, {\em A remark on Fourier transforms}, Proc. Camb.
Phil. Soc. {\bf 32} (1936), 321-327

\bibitem{Zygmund} A.~Zygmund, \it Trigonometric Series, \rm
Cambridge University Press, Cambridge 2002 (3rd edition)




\end{thebibliography}
\end{document}